\documentclass[12pt,reqno]{amsart}
\usepackage{amsmath,amsfonts,amsthm,bbm,amssymb,color,setspace}
\usepackage[usenames,dvipsnames]{xcolor}
\usepackage[T1]{fontenc}
\usepackage{enumerate}
\usepackage{hyperref}
\hypersetup{
     colorlinks   = true,
     citecolor    = blue,
     linkcolor    = blue
}

\usepackage{pdfsync}{\tiny }

\usepackage[left=1in, right=1in, top=1.1in,bottom=1.1in]{geometry}
\numberwithin{equation}{section}

\newcommand{\E}{\mathrm{E}}
\renewcommand{\P}{\mathrm{P}}
\newcommand{\R}{\mathbb{R}}

\newtheorem{theorem}{Theorem}[section]
\newtheorem{lemma}[theorem]{Lemma}

\newtheorem{proposition}[theorem]{Proposition}

\newtheorem{obs}[theorem]{Remark}

\def\1{{\rm l}\hskip -0.21truecm 1}

\begin{document}
\title{Density estimates for jump diffusion processes}
 \date{\today}

\author{Arturo Kohatsu-Higa, Eulalia Nualart and Ngoc Khue Tran}

\address{Arturo Kohatsu-Higa, Department of Mathematical Sciences, Ritsumeikan University, 1-1-1 Nojihigashi, Kusatsu, Shiga, 525-8577, Japan}
\email{khts00@fc.ritsumei.ac.jp}

\address{Eulalia Nualart, Department of Economics and Business, Universitat Pompeu Fabra and Barcelona Graduate School of Economics,  Ram\'on Trias Fargas 25-27, 08005
	Barcelona, Spain.}
\email{eulalia.nualart@upf.edu}

\address{Ngoc Khue Tran, Department of Natural Science Education, Pham Van Dong University, 509 Phan Dinh Phung, Quang Ngai City, Quang Ngai, Vietnam}
\email{tnkhue@pdu.edu.vn}

\thanks{Eulalia Nualart acknowledges support from the Spanish MINECO grant PGC2018-101643-B-I00 and
Ayudas Fundacion BBVA a Equipos de Investigaci\'on Cient\'ifica 2017. Ngoc Khue Tran acknowledges support from the Vietnam Institute for Advanced Study in Mathematics (VIASM) where a part of this work was done during his visit.}


\date{\today}


\begin{abstract}
  We consider a real-valued diffusion process with a linear jump term driven by a Poisson point process and we assume that the jump amplitudes have a centered density with finite moments.
 We show upper and lower estimates for the density of the solution in the case that the jump amplitudes follow a 
 Gaussian or  Laplacian law.  The proof of the lower bound uses a general expression for the  density of the solution in terms of the convolution  of the density of the continuous part and the jump amplitude density.
 The upper bound uses an upper tail estimate in terms of the  jump amplitude distribution and techniques of the Malliavin calculus in order to bound the density by the tails of the solution. We also extend the  lower bounds  to the multidimensional case.
 \end{abstract}

\maketitle

\section{Introduction and main  results}

Consider the following integral equation with jumps 
\begin{equation} \label{e2}
X^x_t=x+\int_0^t \sigma(X^x_s) dB_s+\int_0^t b(X^x_s) ds+ \sum_{i=1}^{\infty}  Y_i \, 1_{T_i \leq t},
\quad t \geq 0,
\end{equation}
where $x \in \R$ and  $(B_t)_{t\geq0}$ is a standard Brownian motion. The jump amplitude sequence $Y=(Y_i)_{i \geq 1}$ is formed with  i.i.d. random variables which have mean zero, finite moments of all orders and probability density function $\varphi$.
The  jump times $(T_i)_{i \geq 1}$ are the arrival times of a Poisson process $(N_t)_{t \geq 0}$ with rate $\lambda>0$.
All sources of randomness are assumed to be mutually independent. 

The coefficients $\sigma, b:\R \rightarrow \R$ are assumed to be twice differentiable with bounded derivatives of all orders. Set $c_1:=\|b\|_\infty$ and $c_2:=\|\sigma\|_\infty$. Moreover, we assume that $\inf_{y \in \R}\vert \sigma(y) \vert \geq \rho>0$ for some constant $\rho>0$. 

Under these conditions it is well-known that there exists a
unique c\`adl\`ag adapted Markov process $X^x=(X^x_t)_{t \geq 0}$ solution to the integral equation \eqref{e2}, which satisfies that for all $T>0$ and $p>1$,
$$\E \bigg[\sup_{t \in [0, T]} \vert X^x_t \vert^p \bigg] < \infty,$$ see for e.g. \cite[Theorem III.2.32]{JS03}.
Moreover, it is also well-known that for all $t>0$, the law of $X^x_t$ has a density with
respect to the Lebesgue measure on $\R$, that we denote by $f_t(x, \cdot)$, see \cite[Theorem 2.5]{Jacod1}
or \cite[Theorem 11.4.4]{NN18}.

In this paper, we are interested in obtaining upper and lower bound estimates for the density $f$. When equation (\ref{e2}) has no jumps,
Gaussian estimates for the density are well-known.
Indeed, if we denote by $Z^x=(Z^x_t)_{t \geq 0}$ the unique solution to the equation 
\begin{equation*}
Z^x_t=x+\int_0^t \sigma(Z^x_s) dB_s+\int_0^t  b(Z^x_s)  ds
\end{equation*}
and by $p_t(x,\cdot)$ its density function, then 
it is well-known that 
 for all $T>0$, there exist constants
$A_T, a_T>1$ such that for all $t \in (0,T]$ and $y \in \R$, 
\begin{equation} \label{Gaussiantype}
\frac{1}{A_T\sqrt{2 \pi  t}} e^{- \frac{a_T\vert y-x \vert^2}{2t}} 
\leq p_t(x,y) \leq \frac{A_T}{\sqrt{2 \pi t}} e^{-\frac{\vert y-x \vert^2}{2 a_T  t}},
\end{equation}
see for e.g. \cite{Kohatsu2, Kusuoka, Sanchez}. However, in the presence of jumps,
less is known about estimates for the density of the solution to (\ref{e2}). Similar estimates for the density function $ f $ as the ones we obtain in this paper are proved in \cite{KusuokaS} (see also the references therein) for a class of infinite activity L\'evy processes. The main motivation to write this paper is the fact that the model (\ref{e2}) appears in some  insurance problems and  their statistical estimation requires in principle properties of their transition densities. The Laplace transform and Karamata-Tuberian theorems are traditionally used in order to obtain asymptotic results for the density, see for e.g. \cite{ACT12} and the references therein. Here we propose a more direct analysis of the density that replaces the machinery of Laplace transforms.

The aim of this paper is to obtain upper and lower bounds for the density $f$ when the jump amplitudes follow the 
Gaussian or Laplace laws. 
\begin{theorem} \label{gaus1}
Assume that $\varphi$ is the centered Gaussian density with variance $\beta>0$. Then for all $T>0$ there exist constants $C_T, c_T>1$ such that for all $t \in (0,T]$ and $x,y \in \R$,
\begin{equation*} \begin{split}
C_T^{-1} \left (e^{-c_{T} \vert y-x \vert \sqrt{\ln_+ (\frac{\vert y-x\vert }{t})}}+\frac{{\bf 1}_{x=y}}{\sqrt{t}}\right ) \leq  f_t(x,y) \leq \frac{C_{T}}{\sqrt{t}}  e^{-c_T^{-1} \vert y-x \vert  \sqrt{\ln_+ (\frac{\vert y-x \vert}{t})}},
\end{split}
\end{equation*}
where $\ln_+(x)=\max(\ln x,0)$. 
\end{theorem}

\begin{theorem} \label{exp1}
Assume that $\varphi$ is the centered Laplace density with scale parameter $1/\mu$ where $\mu>0$, that is, 
\begin{align}
	\label{eq:L}
	\varphi(z)=\frac12 \mu e^{-\mu \vert z \vert}.
\end{align}
Then for all $T>0$ there exist constants $C_T,  c_T>1$ such that for all $t \in (0,T]$ and $x,y \in \R$,
\begin{equation*}
C_T^{-1} \left(e^{-c_T \vert y-x \vert}+ \frac{{\bf 1}_{x=y}}{\sqrt{t}}\right) \leq f_t(x,y)\leq \frac{C_{T}}{\sqrt{t}} e^{-c_T^{-1}\vert y-x \vert}.
\end{equation*}
\end{theorem}

The proof of these two theorems that may be applied to other probability density functions $\varphi$.
Indeed, first, we obtain an expression for the  density $f$ in terms of a convolution  of the density of the continuous part $p$ and the jump amplitude density $\varphi$. This expression is given in Proposition \ref{t1} below and turns out to be very suitable in order to obtain lower bounds for the density.
Second, we show  an upper tail estimate for the solution to equation (\ref{e2}) in terms of the jump distribution of $Y$ that will be crucial for the upper bounds, see Proposition \ref{p1} below. Finally, in Proposition \ref{p2} below we appeal to the techniques of the  Malliavin calculus in order to bound the density $f_t$ in terms of the tail probabilities of $ X_t^x $. These three results are proved for a general jump amplitude density $\varphi$. Then, we will show how to apply these general results for the  two particular cases of $\varphi$ defined in Theorems \ref{gaus1} and \ref{exp1}.

The rest of the paper is organized as follows. Section 2 presents the key preliminary results explained above for a general density $\varphi$. Section 3 is devoted to the proofs of Theorems \ref{gaus1} and \ref{exp1}. Finally, in Section 4 we explore how the main results extend in the multidimensional case.

\section{Preliminary results}

In this section we present some preliminary results that will be crucial for the proof of Theorems \ref{gaus1} and \ref{exp1}.

We start proving an expression for the density that will be  very suitable in order to establish lower bounds.
For any $t>0$ we consider the random variable $Z^x_t+Y$ and we denote by $q_t(x,y)$ its probability density function, where recall that $Z_t^x$ is the continuous part of $X_t^x$ and has density $p$, and $Y$ is the jump amplitude which has density $\varphi$. As $ Y $ and $ Z_t^x $ are independent then we have that $ q_t(x,y)=(p_t*\varphi)(x,y) $ where $ * $ denotes the space convolution of the functions $ p_t(x,\cdot) $ and $ \varphi(\cdot) $. That is,
\begin{align}
\label{q}
q_t(x,y)=\int_{\R}
\varphi(y-v) p_t(x,v) dv.
\end{align}

The following result gives an expression for the density of $X^x_t$ in terms of $q$.
\begin{proposition} \label{t1}
For any $t>0$ and $x,y \in \R$, the density $f_t(x,y)$ of $X^x_t$ solution to equation \textnormal{(\ref{e2})}
can be written as (here $ t_0\equiv 0 $):
\begin{equation} \label{density}
\begin{split}
	& f_t(x,y)=p_t(x,y)e^{-\lambda t}+\sum_{n=1}^{\infty} \int_{t_1< \cdots < t_n<t<t_{n+1}} \!\!\!\!\!\!\!\!\!\!\!\!\!\!\!\! q_{t_1-t_0}* \cdots  * q_{t_n-t_{n-1}}* p_{t-t_n}(x,y)
	\lambda^{n+1} e^{-\lambda t_{n+1}} dt_1 \cdots dt_{n+1}.
\end{split}
\end{equation}
\end{proposition}

\begin{proof}Recall that the probability density function of the random vector $(T_1, T_2,\ldots, T_n)$ has a density function given by
$$
g_{T_1, T_2,\ldots, T_n}(t_1,\ldots, t_n)=\lambda^n e^{-\lambda t_n}1_{0\leq t_1< \cdots < t_{n-1}<t_n}(t_1,\ldots, t_n).
$$
Therefore,
\begin{equation} 
	\label{eq:2}
f_t(x,y)=\E[\delta_y(X^x_t)]=\E[\delta_y(Z^x_t)\, 1_{t<T_1}]+\sum_{n=1}^{\infty} \E[\delta_y(X^x_t) \,1_{T_1< \cdots < T_n<t<T_{n+1}}], 
\end{equation}
where for $n \geq 1$, $X^{x,n}_t\equiv X^{x,n}_t(t_1,...,t_n)$ denotes the solution to the following integral equation with deterministic jump times
\begin{equation*}
X^{x,n}_t=x+\int_0^t \sigma(X^{x,n}_s) dB_s+\int_0^t b(X^{x,n}_s) ds+\sum_{i=1}^{n} Y_i \, 1_{t_i \leq t}. 
\end{equation*}
Remark that an abuse of notation is used here when we write the delta distribution function $\delta_y(x)$. 
The formal argument can be obtained by proper approximation arguments which are left to the reader.

As $ Z^x_t $ and $ T_1 $ are independent, we have $ \E[\delta_y(Z^x_t)\, 1_{t<T_1}] =p_t(x,y) \P(t<T_1)$. Similarly, 
using Chapman-Kolmogorov's equation, we can write 
\begin{equation*}
	\E[\delta_y(X^{x,n}_t)] = q_{t_1-t_0}* \cdots  * q_{t_n-t_{n-1}}* p_{t-t_n}(x,y).
\end{equation*}
Using again the independence between the jump times and the other random components in $ X^{x,n} $, we obtain 
\begin{align*}
	\E[\delta_y(X^x_t) \,1_{T_1< \cdots < T_n<t<T_{n+1}}]=
	\int_{t_1< \cdots < t_n<t<t_{n+1}} \ \E[\delta_y(X^{x,n}_t)]  \lambda^{n+1} e^{-\lambda t_{n+1}} dt_1 \cdots dt_{n+1}.
\end{align*}

Plugging these results into \eqref{eq:2}, the result follows.
\end{proof}

\begin{obs}
Observe that in the linear case, that is, $b=0$ and $\sigma=1$,  \textnormal{(\ref{density})} reads as
\begin{equation*} 
f_t(x,y)=e^{-\lambda t} \sum_{n=0}^\infty (\Phi_t \ast \varphi^{\ast n}) (y-x) \frac{(\lambda t)^n}{n!}
\end{equation*}
where $\Phi_t$ denotes the $N(0,t)$ density.
\end{obs}

The second result of this section is an upper bound for the tail probability of $X_t^x$ in terms of the distribution of the jump amplitude $Y$. This result is an extension of Lemma 26.4 in \cite{Sato}, where a similar tail probability estimate is obtained for L\'evy  processes. Let us first introduce some notation which is similar to that in \cite{Sato}.
Define
\begin{align*}
	C:=\{ u \in \R: \E[e^{ u Y}]<\infty\} \quad \text{ and } \quad  
	s:=\sup (C).
\end{align*}
Note that  $C$ is an interval and assume that $s>0$.
Set, for $u \in C$,
$$
\Psi(u):=\frac12 u^2 c_2^2  +\lambda\E[e^{ u Y}-1].
$$
Then
$$\Psi'(u)=u c_2^2  +\lambda \E[Y e^{u Y }]$$
and
$$\Psi''(u)=c_2^2  +\lambda \E[Y^2 e^{u Y}].
$$
Notice that $\Psi \in C^{\infty}$, and $\Psi''>0$ in the interior of $C$.

Let $u=\theta(\xi)$ be the inverse function of $\xi=\Psi'(u)$, that is,
$$
\xi=\theta(\xi) c_2^2  +\lambda \E[Y e^{\theta(\xi) Y}], \qquad \xi \in (0, \Psi'(s{-})).
$$

\begin{proposition} \label{p1} For all $t>0$ and $x,y \in \R$ such that $\frac{\vert y-x \vert}{t} -c_1\in (0, \Psi'(s{-}))$,
we have
\begin{equation*}
\P ( \vert X^x_t-x\vert  > \vert y-x\vert) \leq 2 e^{-t\int_{0}^{\frac{\vert y-x \vert}{t}-c_1} \theta(\xi) d\xi}.
\end{equation*}
\end{proposition}

\begin{proof}
Let $u \in (0,s)$ and let $(M_t)_{t\geq 0}$ denote the It\^o-L\'evy process given by
\begin{equation*}
\begin{split} 
M_t&=u \int_0^t \sigma(X^x_s) dB_s+u\sum_{i=1}^{\infty} Y_i \, 1_{T_i \leq t} -\frac12  u^ 2 \int_0^t \sigma^2 (X^x_s) ds -\lambda t \int_{\R} 
(e^{uz}-1)\varphi(z) dz.
\end{split}
\end{equation*}
Observe that
$$
u(X_t^x-x)=M_t+u\int_0^t b(X^x_s) ds+\frac12  u^ 2 \int_0^t \sigma^2 (X^x_s) ds+\lambda t\int_{\R} 
(e^{uz}-1)\varphi(z) dz.
$$
By It\^o's formula,
\begin{equation*} 
\begin{split}
e^{M_t}&=1+ u \int_0^t e^{M_s} \sigma(X^x_s) dB_s+\sum_{i=1}^{\infty} e^{M_{T_{i-}}} \left( e^{uY_i}-1\right)1_{T_i \leq t}-
\lambda\int_0^t\int_{\R} e^{M_{s-}} (e^{uz}-1)\varphi(z) dz ds.
\end{split}
\end{equation*}
In particular, $M_t$ is a martingale and $\E[e^{M_t}]=1$. 

Using Markov's inequality and the fact that $\sigma$ and $b$ are bounded, we have that
\begin{equation*}
\begin{split}
\P (  X^x_t-x  > \vert y-x\vert) &=\P ( e^{u (X^x_t-x)}  > e^{u\vert y-x\vert})\\
&=\P ( e^{ M_t+u\int_0^t b(X^x_s) ds+\frac12  u^ 2 \int_0^t \sigma^2 (X^x_s) ds +\lambda\int_0^t\int_{\R} 
		(e^{uz}-1)\varphi(z) dzds } >e^{u\vert y-x\vert})\\
&\leq \P ( e^{ M_t +uc_1t+\frac12  u^ 2 c_2^2 t +\lambda t\int_{\R} 
	(e^{uz}-1)\varphi(z)dz } >e^{u\vert y-x\vert})\\
&= \P \left( e^{ M_t} > e^{(\vert y-x \vert -c_1 t)u -\frac12 u^2 c_2^2 t - \lambda t \E[e^{u Y}-1]  } \right) \\
&\leq \min_{0<u<s} e^{-(\vert y-x \vert -c_1 t)u +t \Psi(u)}=
\min_{0<u<s} e^{t (\Psi(u)-uz)},
\end{split}
\end{equation*}
where $z:=\frac{\vert y-x\vert}{t}-c_1$. The rest of the proof follows as in Lemma 26.4 in \cite{Sato}. Indeed, we have $z\in (0, \Psi'(s_{-}))$. As $\Psi'(u)-z$ changes value from negative to positive at $u=\theta(z)$ and $\Psi(0)=0$, we have
\begin{align*}
\min_{0<u<s}(\Psi(u)-uz)&=\Psi(\theta(z))-\theta(z)z=\int_0^{\theta(z)}\Psi'(u)du-\theta(z)z=\int_{0}^z\xi d\theta(\xi)-\theta(z)z\\
&=z\theta(z)-\lim_{\xi \downarrow 0}\xi\theta(\xi)-\int_{0}^z\theta(\xi)d\xi-\theta(z)z=-\int_{0}^z\theta(\xi)d\xi,
\end{align*}
since 
\begin{align*}
\lim_{\xi \downarrow 0}\xi\theta(\xi)=\lim_{u \downarrow 0}\Psi'(u)u=0.
\end{align*}
This shows the upper bound for $\P (  X^x_t-x  > \vert y-x\vert)$.

Proceeding exactly along the same lines, we can consider the martingale $-M_t$ and show the same upper bound for $\P ( -( X^x_t-x)  > \vert y-x\vert)$.
Thus, the desired result follows.
\end{proof}

The last result of this section is an upper bound for the density of $X_t^x$ in terms of its upper tails.
For this, we appeal to techniques of the Malliavin calculus with respect to the Brownian motion $B$.
We denote by $D$ the Malliavin derivative operator with respect to $B$ and by $\mathbb{D}^{2, \infty}$ the Sobolev space of twice differentiable random variables with finite moments of all orders. See the monographs \cite{N06} or \cite{NN18} for the precise definitions.
The next result is classical and shows that the solution to \eqref{e2} is twice differentiable in the Malliavin sense and gives an expression for the Malliavin derivatives, see \cite{Jacod1, Jacod2} and  \cite[Theorem 11.4.3]{NN18} for its proof.
\begin{lemma} \label{mal}
For all $t>0$ and $x \in \R$, $X^x_t$ belongs to $\mathbb{D}^{2,\infty}$ and the Malliavin derivative $(D_rX^x_t, r\leq t)$ satisfies the following linear equation
\begin{equation*}
\begin{split}
D_rX^x_t=\sigma(X^x_r)+\int_r^t \sigma'(X^x_s)D_rX^x_s dB_s+\int_r^t b'(X^x_s)D_rX^x_s ds,
\end{split}
\end{equation*}
for $r\leq t$, a.e., and $D_rX^x_t=0$ for $r>t$, a.e.
Moreover, for all $p>1$,
\begin{equation*}
\sup_{r\in[0,T]}\E\left[\sup_{t\in [r,T]}\left\vert D_rX^x_t\right\vert^p\right]<\infty.
\end{equation*} 
Furthermore, the iterated  Malliavin derivative  $(D^2_{r_1,r_2}{X^x_t}, r_1\vee r_2 \leq t)$, satisfies the equation
\begin{equation*} 
\begin{split}
D^2_{r_1,r_2}X^x_t&=D_{r_1}\sigma(X^x_{r_2})+D_{r_2}\sigma(X^x_{r_1})+\int_{r_1\vee r_2}^tD^2_{r_1,r_2}(\sigma(X^x_s))\, dB_s\\
&\quad+\int_{r_1\vee r_2}^tD^2_{r_1,r_2}(b(X^x_s))\, ds,
\end{split}
\end{equation*}
for $r_1\vee r_2\leq t$, a.e., and $D^2_{r_1,r_2}X^x_t=0$ otherwise.
Moreover, for all $p>1$,
\begin{equation*}
\sup_{r_1,r_2\in[0,T]}\E\left[\sup_{r_1\vee r_2\leq t\leq T}\left\vert D^2_{r_1,r_2}X^x_t\right\vert^p\right]<\infty.
\end{equation*} 
\end{lemma}

We are now ready to state the last result of this section  which bounds the density $f_t(x,y)$ of $X_t^x$ in terms of its upper tail probability. Since the jump  term is linear, the proof follows exactly along the same lines as for a continuous diffusion (see \cite{N06}). For the sake of completeness, we provide the proof.
\begin{proposition} \label{p2}
For any $T>0$ and $q>1$, there exists a constant $C_{q,T}>0$ such that for all $t \in (0,T]$ and $x,y\in\R$,
\begin{equation*}
 f_t(x,y) \leq  \frac{C_{q,T}}{\sqrt{t}} (\P ( \vert X^x_t-x\vert  > \vert y-x\vert))^{1/q}.
\end{equation*}
\end{proposition}

\begin{proof}
Appealing to Proposition 2.1.2 in \cite{N06}, we have
\begin{equation} \label{eqf}
\begin{split}
f_t(x,y)&\leq c_{q,\alpha,\beta}\left(\P ( \vert X^x_t-x\vert  > \vert y-x\vert)\right)^{1/q}\bigg(\E\left[\Vert DX^x_t \Vert_H^{-1}\right]+\left(\E\left[\Vert D^2X^x_t \Vert_{H\otimes H}^{\alpha}\right]\right)^{1/\alpha}\\
&\qquad\times\left(\E\left[\Vert DX^x_t \Vert_{H}^{-2\beta}\right]\right)^{1/\beta}\bigg),
\end{split}
\end{equation}
for some constant $c_{q,\alpha,\beta}>0$, where  $\frac{1}{q}+\frac{1}{\alpha}+\frac{1}{\beta}=1$ and $H=L^2([0,t];\R)$.

Using Lemma \ref{mal} and H\"older's inequality, it is easy to check that for any $\alpha>1$,
\begin{equation} \label{eqf1}
\begin{split}
\left(\E\left[\Vert D^2X^x_t \Vert_{H\otimes H}^{\alpha}\right]\right)^{1/\alpha}\leq C_{\alpha,T}t \leq C'_{\alpha,T}\sqrt{t}.
\end{split}
\end{equation}

Now, we denote by $J_t^x:=\partial_x X_t^x$ the derivative of $X_t^x$ with respect to the initial condition $x$, which satisfies the linear equation
\begin{equation*}
\begin{split}
J_t^x=1+\int_0^t \sigma'(X^x_s)J_s^x dB_s+\int_0^t b'(X^x_s)J_s^x ds.
\end{split}
\end{equation*}
By It\^o's formula, the inverse $H^x_t:=(J_t^x)^{-1}$ satisfies the linear equation
\begin{equation*}
\begin{split}
H^x_t=1-\int_0^t \sigma'(X_s^x)H^x_s dB_s -\int_0^t(b'(X_s^x)-(\sigma'(X_s^x))^2 ) H^x_sds
\end{split}
\end{equation*}
Using the assumptions on $b$ and $\sigma$ and Gronwall type arguments, we have  that for all $p\geq 1$,
\begin{equation*}
\begin{split}
\E\Big[\sup_{t\in [0,T]} \vert J_t^x \vert^p\Big]\leq C_{p,T},\qquad \E\Big[\sup_{t\in [0,T]} \vert  H_t^x \vert^{p}\Big]\leq C_{p,T}.
\end{split}
\end{equation*}
Moreover, we have  that
\begin{equation*}
\begin{split}
D_r X_t^x =  J_t^x H_r^x \sigma(X_r^x) 1_{[0,t]}(r).
\end{split}
\end{equation*}
Consequently, for any $p \geq 1$,
\begin{equation} \label{eqf2}
\begin{split}
\E\left[\Vert DX_t \Vert_H^{-p}\right]&=\E\left[\left(\int_0^t(J_t^x)^2(H_r^x)^{2} \sigma^2(X_r^x)dr\right)^{-p/2}\right]\\
&\leq \dfrac{1}{\rho^{p}t^{p/2}} \E\left[\sup_{0\leq r\leq T}\vert J_r^x\vert ^{p} \sup_{0\leq r\leq T} \vert H_r^x \vert^{p}\right]\leq \dfrac{C_{p,T}}{\rho^{p}t^{p/2}}.
\end{split}
\end{equation}
Thus, substituting (\ref{eqf1}) and (\ref{eqf2}) into (\ref{eqf}), we conclude the proof of the upper bound.
Observe that Lemma \ref{mal}, (\ref{eqf2}) and  criterion in \cite[Theorem 2.1.1]{N06} imply the existence of the density.
\end{proof}

\section{Proof of Theorems \ref{gaus1} and \ref{exp1}}

\begin{proof}[Proof of Theorem \ref{gaus1}]
Assume that $\varphi$ is $N(0, \beta)$. We start proving the upper bound. By Proposition \ref{p2}, it suffices to apply Proposition \ref{p1}. For $u>0$, we have that
\begin{align*}
\Psi(u)=\frac12 u^2 c_2^2  +\lambda (e^{\frac{u^2 \beta }{2}}-1).
\end{align*}
Thus, $\Psi'(u)=u (c_2^2  +\lambda \beta {e^{\frac{u^2 \beta}{2}}})$ and its inverse function $u=\theta(\xi)$ satisfies
$$
\xi=\theta(\xi) \big(c_2^2  +\lambda \beta  e^{\frac{\theta(\xi)^2 \beta }{2}}\big), \quad \xi \in (0, \infty).
$$ 

Let us now estimate $\theta(\xi)$. Observe that $\theta(\xi)\uparrow \infty$ as $\xi\uparrow\infty$ and for $\alpha < \frac{2}{\beta }$, $\xi e^{-\frac{\theta^2(\xi)}{\alpha}}\to 0$ as $\xi\uparrow\infty$.
Hence, there is $\xi_1 >1$ such that for all $\xi >\xi_1$, $\xi e^{-\frac{\theta^2(\xi)}{\alpha}}<\frac12$ and
$c_1 e^{-\frac{\theta^2(\xi)}{\alpha}}<\frac12$. Thus, for $\xi>\xi_1$, $(\xi+c_1) e^{-\frac{\theta^2(\xi)}{\alpha}}<1$ and $\theta(\xi) > \sqrt{\alpha \ln (\xi+c_1)}$. Therefore, by Proposition \ref{p1}, for $\frac {\vert y-x \vert}t-c_1> \xi_1$,
\begin{equation*}
\P (  \vert X^x_t-x\vert  >\vert y-x\vert) \leq 2 e^{-t\int_{\xi_1}^{\frac{\vert y-x \vert}t-c_1} \sqrt{\alpha \ln (\xi+c_1)} d\xi}=2 e^{-t\int_{\xi_1+c_1}^{\frac{\vert y-x \vert}t} \sqrt{\alpha \ln (\xi)} d\xi}. 
\end{equation*}
As for $z$ sufficiently large $\int_{\xi_1+c_1}^{z} \sqrt{\ln \xi} d\xi > \frac z2 \sqrt{\ln z}$, the desired result follows.
In the case that $0<\frac{\vert y-x\vert }t-c_1\le \xi_1$ then the result follows trivially.
In the case that $\frac{\vert y-x \vert}t-c_1<0$ then the probability of the reverse inequality can be bounded in a similar fashion. 

We next prove the lower bound. 
Using \eqref{q}  and the  lower bound in \eqref{Gaussiantype}, we obtain that for all $t \in (0,T]$ and $x,y \in \R$, 
\begin{equation*}
\begin{split}
q_t(x,y) &\geq 
%\int_{\R} \frac{1}{\sqrt{2\pi \beta}} e^{-\frac{\vert y-v \vert^2}{2\beta}}\frac{1}{A_T\sqrt{2\pi t}} e^{-\frac{a_T\vert x-v \vert^2}{2t}} \, dv \\
%&=
\frac{1}{A_T\sqrt{a_T 2\pi (a_T^{-1}t+\beta )} } e^{-\frac{\vert y-x \vert^2}{2(a_T^{-1}t+  \beta )}}.
\end{split}
\end{equation*}

Therefore, using (\ref{density}), we obtain for $r:=\vert y-x \vert$ and  $C_T=\frac{1}{A_T\sqrt{a_T}}$
\begin{align} \label{au2}
f_t(x,y)\geq e^{-\lambda t}\sum_{n=0}^{\infty} \left(C_T\right)^{n}\frac{C_T}{\sqrt{ 2\pi (a_T^{-1}t+n\beta )} } e^{-\frac{r^2}{2(a_T^{-1}t+n \beta )}}\frac{(\lambda t)^n}{n!}.
\end{align}

Observe that it suffices  to assume that $\frac{r}{t} \geq e$, otherwise the  bound follows trivially, since taking $n=0$ yields to
$$
f_t(x,y)\geq e^{-\lambda t} \frac{1}{A_T\sqrt{2\pi t} } e^{-\frac{r^2}{2a_T^{-1}t}} \geq \frac{1}{A_T\sqrt{2\pi T}} e^{-\lambda  T}e^{-\frac{T e^2}{2a_T^{-1}}},
$$
from which the desired lower bound follows for $r < t e$. Observe also that if $r=0$, we get the lower bound
\begin{equation} \label{op}
f_t(x,y)\geq e^{-\lambda T} \frac{1}{A_T\sqrt{2\pi t} }.
\end{equation}

By Stirling's formula, there exists a constant $K > 1$ such that for all $n > K$, it holds that
\begin{equation*}
\left| \frac{\sqrt{2\pi} e^{n \ln(n) -n+\frac12\ln(n)}}{n!} -1\right| < \frac12.
\end{equation*}
This implies that for all $n>K$,
$$
\frac{(\lambda t)^n}{n!}>\frac 1{2\sqrt{2\pi}} e^{-n \ln\left(\frac{n}{\lambda t}\right) +n-\frac12\ln(n)}>\frac 1{2\sqrt{2\pi}} e^{-n \ln\left(\frac{n}{\lambda t}\right) -\frac12\ln(n)}.
$$
Then, substituting this into (\ref{au2}), we get that
\begin{align} \nonumber
f_t(x,y) &\geq \frac{C_T}{4 \pi} e^{-\lambda t}\sum_{n>K} e^{n\ln(C_T)-\frac{1}{2} \ln(a_T^{-1}t+n \beta)-\frac{r^2}{2(a_T^{-1} t+n \beta)}-n \ln\left(\frac{n}{\lambda t}\right)-\frac12\ln(n)} \\
& \geq\frac{C_T}{4 \pi} e^{-\lambda t} \sum_{n>K} e^{-n \vert \ln(C_T)\vert-\frac{1}{2} \ln(a^{-1}_T T+n\beta)-\frac{r^2}{2n\beta}-n \ln\left(\frac{n}{\lambda t}\right)-\frac12\ln(n)}.\label{a4}
\end{align}
We next consider two different cases  according to $\frac{r}{\sqrt{\ln (r/t)}}>K+1$
and $\frac{r}{\sqrt{\ln (r/t)}}\le K+1$. 

In the first case, we set $n=[\frac{r}{\sqrt{\ln (r/t)}}]>K>1$.
Note that 
$x-1\leq [x]\leq x$ for any $ x\in\mathbb{R} $, then 
%$$\frac{r}{\sqrt{\ln (r/t)}}-1 \leq [\frac{r}{\sqrt{\ln (r/t)}}] \leq \frac{r}{\sqrt{\ln (r/t)}}.$$ Hence, 
using (\ref{a4}), we obtain
\begin{equation*} \begin{split}
 &f_t(x,y)\\
&\geq \frac{C_T}{4 \pi} e^{-\lambda t} e^{- \frac{r}{\sqrt{\ln (r/t)}}\vert \ln(C_T)\vert-\frac{1}{2} 
\ln\left(c'_T\frac{r}{\sqrt{\ln (r/t)}}\right)-\frac{r^2}{2\left(\frac{r}{\sqrt{\ln (r/t)}}-1\right)\beta }-
\frac{r}{\sqrt{\ln (r/t)}} \ln\left(\frac{r}{\lambda t\sqrt{\ln (r/t)}}\right)-\frac12\ln \left(\frac{r}{\sqrt{\ln (r/t)}}\right)} \\
&\geq \frac{C_T}{4 \pi} e^{-\lambda T} e^{-c_T r \sqrt{\ln (r/t)}},
\end{split}
\end{equation*}
for some constant $c_T>0$,
where $c_T'=a^{-1}_T T+\beta$. In the last inequality we have used the following inequalities
\begin{enumerate}
	\item $\displaystyle{ \frac{r^2}{\frac{r}{\sqrt{\ln (r/t)}}-1}\leq \frac{K+1}{K} r \sqrt{\ln (r/t)},} $
	\item $ \displaystyle{\frac{r}{\sqrt{\ln (r/t)}}\leq r \leq r\sqrt{\ln (r/t)},\quad \text{since } \frac{r}{t}\geq e,} $
	\item$ \displaystyle{\frac{r}{\sqrt{\ln (r/t)}}\ln\left(\frac{r}{\lambda t\sqrt{\ln (r/t)}}\right) \leq \frac{r}{\sqrt{\ln (r/t)}} \ln\left(\frac{r/t}{\lambda }\right)\leq r\sqrt{\ln (r/t)}(1+\vert \ln(\lambda) \vert).} $
\end{enumerate}

For the second case $\frac{r}{\sqrt{\ln (r/t)}}\le K+1$ the conclusion follows easily since taking $n=0$ in (\ref{au2}) yields 
$$
f_t(x,y)\geq e^{-\lambda t} \frac{1}{A_T\sqrt{2\pi T} } e^{-\frac{r^2}{2a_T^{-1}T}}\geq 
e^{-\lambda T} \frac{1}{A_T\sqrt{2\pi T} } e^{-(K+1)\frac{r \sqrt{\ln (r/t)}}{2a_T^{-1}T}}.
$$
The proof is now completed.
\end{proof}

\begin{proof}[Proof of Theorem \ref{exp1}]
We first prove the upper bound. As above it suffices to apply Propositions \ref{p1} and \ref{p2}. A direct calculation using  \eqref{eq:L} gives that for all $-{\mu}<u<{\mu}$, 
$$
\Psi(u)=\frac12 u^2 c_2^2  +\lambda \left(\frac{\mu}{2}\left(\frac{1}{u+\mu}-\dfrac{1}{u-\mu}\right) -1\right).
$$
Thus,
$$
\xi=\theta(\xi) c_2^2+ \frac{\lambda \mu }{2} \left(\frac{1}{(\theta(\xi) -\mu)^2}-\frac{1}{(\theta(\xi) +\mu)^2}\right),
\quad \xi \in (0,\infty).
$$

Let us now estimate $\theta(\xi)$. Observe that $\theta(\xi)\uparrow{\mu}$ as $\xi\uparrow\infty$. 
This implies that $\xi( (\theta(\xi)-\mu)^2 \wedge (\theta(\xi)+\mu)^2)$ 
converges to $\frac{\lambda \mu }{2}$ as $\xi\uparrow \infty$.
Hence, there exists $\xi_1>0$ such that for all $\xi >\xi_1$, $\sqrt{\xi}({\mu}-\theta(\xi))-{\sqrt{\lambda \mu/2}}<1$, and thus, for all $\xi >\xi_1$,
$\theta(\xi)>{\mu}-\frac{1+\sqrt{\lambda \mu/2}}{\sqrt{\xi}}$. Therefore, by Proposition \ref{p1}, for $\frac {\vert y-x \vert}{t}-c_1> \xi_1$,
\begin{equation*}
\P (\vert X^x_t-x\vert>\vert y-x\vert) \leq 2 e^{-t\int_{\xi_1}^{\frac {\vert y-x \vert}t-c_1} \left(\mu-\frac{1+\sqrt{\lambda \mu/2}}{\sqrt{\xi}}\right) d\xi}. 
\end{equation*}
Choosing  $\frac{1} {\sqrt{\xi_1}} >\frac{\mu}{2(1+\sqrt{\lambda \mu/2})}$, gives

\begin{align*}
\int_{\xi_1}^{\frac {\vert y-x \vert}t-c_1} (\frac{\mu}{1+\sqrt{\lambda \mu/2}}-\frac{1}{\sqrt{\xi}}) d\xi >\frac{\mu}{2(1+\sqrt{\lambda \mu/2})}(\frac {\vert y-x \vert}t-c_1-\xi_1).
\end{align*}  
In the case that $\frac{\vert y-x \vert}t-c_1\le \xi_1$ then the result follows trivially. This finishes the proof of the upper bound.

Next, we prove the lower bound. Using \eqref{q} and the lower bound in \eqref{Gaussiantype}, we get that for all $t \in (0,T]$ and $x,y \in \R$, 
\begin{equation*} 
\begin{split}
q_t(x,y) &\geq \frac{\mu}{2 }\int_{\R} 
e^{-\mu \vert z \vert} \frac{1}{A_T\sqrt{2\pi t}}e^{-\frac{\vert z-(y-x)\vert^2}{2a_T^{-1} t}} dz \\
&=\frac{\mu}{2 A_T\sqrt{2\pi t}}e^{\frac{\mu^2 ta_T^{-1}}{2}} \bigg(e^{-(y-x) \mu}\int_{0}^{\infty} 
e^{-\frac{\vert z-(y-x-\mu t a_T^{-1})\vert^2}{2a_T^{-1} t}}dz  +e^{(y-x) \mu} \int_{-\infty}^{0} 
e^{-\frac{\vert z-(y-x+\mu t a_{T}^{-1})\vert^2}{2a_T^{-1} t}}dz\bigg) \\
&\geq C_T e^{\frac{\mu^2 ta_T^{-1}}{2}} \bigg(e^{-(y-x) \mu} \left({\bf 1}_{y-x-\mu t a_T^{-1}<0} \frac12 e^{-\frac{\vert y-x-\mu t a_T^{-1}\vert^2}{a_T^{-1} t}}+\frac12 {\bf 1}_{y-x-\mu t a_T^{-1}\geq  0 }\right) \\
&\qquad \qquad  +e^{(y-x) \mu} \left({\bf 1}_{y-x+\mu t a_T^{-1}\geq 0} \frac12  e^{-\frac{\vert y-x+\mu t a_T^{-1}\vert^2}{a_T^{-1} t}}+\frac12 {\bf 1}_{y-x+\mu t a_T^{-1}< 0 } \right) \bigg),
\end{split}
\end{equation*}
where $C_T=\frac{\mu}{2A_T \sqrt{2a_T}}$. Observe that in order to get the last inequality when $z>0$, we have used the fact that the integral of a Gaussian density with a non-negative mean on the positive axis is lower bounded by $\frac12$. On the other hand, in the case that the mean is negative, we have used the inequality $\vert a-b\vert^2\leq 2 (\vert a\vert^2+\vert b \vert^2)$, valid for all $a,b \in \R$ and the fact that the integral of a Gaussian density with zero mean on the positive axis equals $\frac12$. We have applied a similar argument for the case $z<0$.

Expanding the square appearing in the two exponentials, yields 
\begin{equation*} 
\begin{split}
q_t(x,y)
&\geq \frac{C_T}{2}e^{-\frac{\mu^2 ta_T^{-1}}{2}} \bigg(e^{-\frac{\vert y-x\vert^2}{a_T^{-1} t}} 
\left(e^{(y-x) \mu}{\bf 1}_{y-x<\mu t a_T^{-1}}+e^{-(y-x) \mu}{\bf 1}_{y-x\geq  -\mu t a_T^{-1}}\right)\\
&\qquad \qquad+e^{-(y-x) \mu} {\bf 1}_{y-x\geq \mu t a_T^{-1}}
 +e^{(y-x) \mu} {\bf 1}_{y-x<-\mu t a_T^{-1}}\bigg) \\
 &\geq \frac{C_T}{2}e^{-\frac{\mu^2 ta_T^{-1}}{2}} \bigg(e^{-\frac{\vert y-x\vert^2}{a_T^{-1} t}} 
\bigg(e^{-\vert y-x\vert \mu}{\bf 1}_{\vert y-x\vert  \leq\mu t a_T^{-1}}\\
&\qquad \qquad+e^{-\vert y-x \vert \mu} {\bf 1}_{y-x> \mu t a_T^{-1}}
 +e^{-\vert y-x \vert \mu} {\bf 1}_{y-x<-\mu t a_T^{-1}}\bigg)\bigg) \\
%&\geq \frac{C_T}{2} e^{-\frac{\mu^2 ta_T^{-1}}{2}}\bigg(e^{-2\vert y-x\vert \mu} {\bf 1}_{-\mu t a_T^{-1} \leq   
%y-x <\mu t a_T^{-1}}+e^{-\vert  y-x \vert  \mu} {\bf 1}_{y-x\geq \mu t a_T^{-1}}\\
%&\qquad \qquad +e^{-\vert y-x \vert \mu}{\bf 1}_{y-x<-\mu t a_T^{-1}}\bigg) 
%\\
&\geq \frac{C_T}{2}  e^{-\frac{\mu^2 ta_T^{-1}}{2}}e^{-2\vert y-x\vert \mu}.
\end{split}
\end{equation*}
This implies, using the triangular inequality, that
\begin{equation*} \begin{split}
 q_{t_1}* q_{t_2-t_1}(x, y_2) 
		&\geq \frac{C_T^2}{4}\int_{\R} e^{-\frac{\mu^2 t_1 a_T^{-1}}{2}} e^{-2\vert y_1-x\vert \mu}e^{-\frac{\mu^2 (t_2-t_1) a_T^{-1}}{2}}e^{-2\vert y_2-y_1\vert \mu} dy_1\\
	&=\frac{C_T^2}{4} e^{-\frac{\mu^2 t_2 a_T^{-1}}{2}}\int_{\R} e^{-2\vert y_1\vert \mu}e^{-2\vert y_2-x-y_1\vert \mu} dy_1\\
	&\geq 
%	\frac{C_T^2}{4} e^{-\frac{\mu^2 t_2 a_T^{-1}}{2}}e^{-2\vert y_2-x\vert \mu}\int_{\R} e^{-4\vert y_1\vert \mu} dy_1 \\
%	&=
	{2 \mu C_T^2} e^{-\frac{\mu^2 t_2 a_T^{-1}}{2}}e^{-2\vert y_2-x\vert \mu}.
	\end{split}
	\end{equation*}

Therefore, iterating the above computation, we obtain that for $n \geq 1$,
\begin{equation*} \begin{split}
 q_{t_1}* q_{t_2-t_1}*\cdots *q_{t_n-t_{n-1}}* p_{t-t_n}(x,y) 
	&\geq \tilde{C}_T^n e^{-\frac{\mu^2 t_n a_T^{-1}}{2}}  \int_{\R} e^{-2\vert y_n-x\vert \mu} p_{t-t_n}(y_n,y) dy_n\\
	&\geq \overline{C}_T \tilde{C}_T^{n} e^{-c_Tt} e^{-4 \vert y-x \vert \mu},
	\end{split}
	\end{equation*}
	for some constants $\overline{C}_T, \tilde{C}_T, c_T>0$. Finally, appealing to formula (\ref{density}) yields to
	\begin{align*} 
f_t(x,y)&\geq \overline{C}_T e^{-c_Tt} e^{-4 \vert y-x \vert \mu} e^{-\lambda t}\sum_{n=0}^{\infty} \tilde{C}_T^{n}\frac{(\lambda t)^n}{n!} \geq \overline{C}_T e^{-(c_T+\lambda) T} e^{-4 \vert y-x \vert \mu},
\end{align*}
which proves the lower bound for $x \neq y$.
When $x=y$ it suffices to use formula (\ref{density}) for $n=0$ to obtain the same lower bound as in (\ref{op}).
This completes the proof.
\end{proof}

\section{Extension to the multidimensional case}

The aim of this section is to explain how the results obtained above extend to the multidimensional case. 
The multidimensional version of equation (\ref{e2}) writes as follows:
\begin{equation} \label{e2d}
X^x_t=x+\int_0^t \sigma(X^x_s) dB_s+\int_0^t b(X^x_s) ds+\sum_{i=1}^{\infty} Y_i \, 1_{T_i \leq t},
\quad t \geq 0,
\end{equation}
where $x \in \R^d$, $(B_t)_{t\geq0}$ is a $d$-dimensional standard Brownian motion, $Y=(Y_i)_{i \geq 1}$ is a sequence of $d$-dimensional  i.i.d. random variables  
which have mean zero, finite moments of all orders and probability density function $\varphi$.
The function $b:\R^d \rightarrow \R^d$ and the matrix $\sigma: \R^d \rightarrow \mathcal{M}_{d \times d}$ are $\mathcal{C}^{\infty}$, bounded with bounded partial derivatives of all orders. We also assume that the matrix $\sigma$ is uniformly elliptic,  that is, there exists $ \rho>0 $ such that
$$
\inf_{\xi \in \R^d: \vert  \xi \vert=1} \vert \sigma(y)\xi \vert^2 \geq \rho >0.
$$

Under theses conditions it is well-known that there exists a unique c\`adl\`ag adapted Markov process $X^x=(X^x_t)_{t \geq 0}$ solution to the integral equation \eqref{e2d}, see \cite{JS03}. Moreover, for all $t>0$ the random vector $X_t^x$ possesses a density $f_t(x,\cdot)$ with
respect to the Lebesgue measure on $\R^d$, see \cite{Jacod1}.

As before we denote by $p_t(x,\cdot)$ the density function of the solution to equation \eqref{e2d} with $Y\equiv 0$. Then,
it is well-known that 
 for all $T>0$, there exist constants
$A_T, a_T>1$ such that for all $t \in (0,T]$ and $x,y \in \R^d$, 
\begin{equation} \label{g1}
\frac{1}{A_T(2 \pi  t)^{d/2}} e^{- \frac{a_T\vert y-x \vert^2}{2t}} 
\leq p_t(x,y) \leq \frac{A_T}{(2 \pi t)^{d/2}} e^{-\frac{\vert y-x \vert^2}{2 a_T  t}}.
\end{equation}

The expression for the density obtained in Proposition \ref{t1} can be easily extended in this multidimensional setting as follows. The proof follows exactly as in the one-dimensional case.
\begin{proposition} \label{t1d}
For any $t>0$ and $x,y \in \R^d$, 
\begin{equation*} 
\begin{split}
	&f_t(x,y)=p_t(x,y)e^{-\lambda t}+\sum_{n=1}^{\infty} \int_{t_1< \cdots < t_n<t<t_{n+1}} \!\!\!\!\!\!\!\!\!\!\!\!\!\!\!\! q_{t_1-t_0}* \cdots  * q_{t_n-t_{n-1}}* p_{t-t_n}(x,y)
  \lambda^{n+1} e^{-\lambda t_{n+1}} dt_1 \cdots dt_{n+1}.
\end{split}
\end{equation*}
\end{proposition}

The lower bounds of Theorems \ref{gaus1} and \ref{exp1} also extended to equation (\ref{e2d}) as follows.
\begin{theorem} \label{gaus1b}
Assume that $\varphi$ is the centered $d$-dimensional Gaussian density with covariance matrix $\Sigma$. Then for all $T>0$ there exist constants $C_T, c_T>1$ such that for all $t \in (0,T]$ and $x,y \in \R^d$,
\begin{equation*} 
f_t(x,y) \geq C_T^{-1}\left (e^{-c_{T} \vert y-x \vert \sqrt{\ln_+ (\frac{\vert y-x\vert }{t})}}+\frac{{\bf 1}_{x=y}}{t^{d/2}}\right ).
\end{equation*}
\end{theorem}

\begin{proof}
Using convolution properties for Gaussian densities with the lower bound in (\ref{g1}) yields, for all $t \in (0,T]$ and $x,y \in \R^d$,
\begin{equation*}
\begin{split}
q_t(x,y) &\geq 
%\int_{\R^d} \frac{1}{ (2\pi)^{d/2} \sqrt{\vert \det(\Sigma) \vert }} 
%e^{- \frac{(y-v)^T  \Sigma ^{-1}
% (y-v)}{2}}\frac{1}{A_T a_T^{d/2}(2\pi a_T^{-1} t)^{d/2}} e^{-\frac{\vert x-v \vert^2}{2 a_T^{-1} t}} \, dv \\
%&=
\frac{1}{A_T  (2 \pi a_T)^{d/2} \sqrt{\vert  \det(\Sigma+a_T^{-1} t I) \vert } } e^{- \frac{(y-x)^T (\Sigma+a_T^{-1} t I)^{-1}
 (y-x)}{2}},
\end{split}
\end{equation*}
where $I$ denotes the identity matrix of order $d \times d$.
Therefore, using Proposition \ref{t1d}, we obtain in a similar manner that
\begin{equation*} \begin{split}
f_t(x,y)
&\geq e^{-\lambda t}\sum_{n=0}^{\infty} C_T^n\frac{C_T}{(2 \pi)^{d/2} \sqrt{\vert  \det(n\Sigma+a_T^{-1} t I) \vert } } e^{- \frac{(y-x)^T (n\Sigma+a_T^{-1} t I)^{-1}
 (y-x)}{2}}\frac{(\lambda t)^n}{n!} \\
 &\geq e^{-\lambda t}\sum_{n=1}^{\infty} C_T^n\frac{C_T}{(2 \pi)^{d/2} n^{d/2} \sqrt{\vert  \det(\Sigma+a_T^{-1} T I) \vert } } e^{-\frac{r^2 \Vert (\Sigma)^{-1} \Vert}{2n}} \frac{(\lambda t)^n}{n!},
\end{split}
\end{equation*}
where $r:=\vert y-x \vert$, $C_T=\frac{1}{A_T(a_T)^{d/2}}$, and $\Vert (\Sigma)^{-1} \Vert=\sup_{z\in \R^d: z \neq 0} \frac{\vert (\Sigma)^{-1}z\vert  }{\vert z \vert}$. 

The rest of the proof follows exactly along the  same lines as in the one-dimensional case.
\end{proof}

\begin{theorem} \label{exp1b}
Assume that $\varphi$ is the multivariate centered Laplace density with $\mu=(\mu_1,\ldots,\mu_d)$, $\mu_i>0$, given by
$$
\varphi(z)=\prod_{i=1}^d\frac12 \mu_i e^{-\mu_i \vert z_i \vert}, \quad z=(z_1,\ldots,z_d).
$$
Then for all $T>0$ there exist constants $C_T, c_T>1$ such  that for all $t \in (0,T]$ and $x,y \in \R^d$,
\begin{equation*}
f_t(x,y)\geq C_T^{-1} \left(e^{-c_T \vert y-x \vert}+ \frac{{\bf 1}_{x=y}}{t^{d/2}}\right).
\end{equation*}
\end{theorem}

\begin{proof} 
We start proving the result for $d=2$. We proceed similarly as in the case $d=1$, but in this case  we need to split the integral into the 4 quadrants of $\R^2$.
Using (\ref{q}) and the lower bound in (\ref{g1}), we get that for all $t \in (0,T]$ and $x,y \in \R^2$, with $ \hat{\mu}:=(\mu_1,-\mu_2) $
\begin{align}
&q_t(x,y) 
\geq \frac{\mu_1 \mu_2}{8\pi tA_T }\int_{\R^2} 
e^{-\mu_1 \vert z_1 \vert} e^{-\mu_2 \vert z_2 \vert} e^{-\frac{\vert z-(y-x)\vert^2}{2a_T^{-1} t}} dz \nonumber\\
&= \frac{\mu_1 \mu_2}{8\pi tA_T }e^{\frac{a_T^{-1} t\vert  \mu \vert^2}{2 }} \bigg(e^{- \langle y-x, \mu \rangle}\int_{\mathbb{R}_+^2} dz
e^{-\frac{\vert z-(y-x-a_T^{-1} t \mu)\vert^2}{2a_T^{-1} t}} +e^{- \langle y-x, \hat \mu \rangle} \int_{\mathbb{R}_+\times \mathbb{R}_-}dz
e^{-\frac{\vert z-(y-x-a_T^{-1} t \hat \mu)\vert^2}{2a_T^{-1} t}}\nonumber
\\
&\qquad +e^{\langle y-x, \hat{\mu} \rangle}\int_{\mathbb{R}_-\times \mathbb{R}_+}dz
e^{-\frac{\vert z-(y-x+a_T^{-1} t \hat \mu)\vert^2}{2a_T^{-1} t}}+e^{\langle y-x,\mu \rangle} \int_{\mathbb{R}_-^2}dz
e^{-\frac{\vert z-(y-x+a_T^{-1}t \mu)\vert^2}{2a_T^{-1} t}}\bigg).
\label{eq:8}
\end{align}
%\begin{equation} 
%\begin{split}
%&q_t(x,y) 
%\geq \frac{\mu_1 \mu_2}{8\pi tA_T }\int_{\R^2} 
%e^{-\mu_1 \vert z_1 \vert} e^{-\mu_2 \vert z_2 \vert} e^{-\frac{\vert z-(y-x)\vert^2}{2a_T^{-1} t}} dz \\
%&= \frac{\mu_1 \mu_2}{8\pi tA_T }e^{\frac{a_T^{-1} t\vert  \mu \vert^2}{2 }} \bigg(e^{- \langle y-x, \mu \rangle}\int_{\mathbb{R}_+^2} dz
%e^{-\frac{\vert z-(y-x-a_T^{-1} t \mu)\vert^2}{2a_T^{-1} t}} +e^{- \langle y-x, \hat \mu \rangle} \int_{\mathbb{R}_+\times \mathbb{R}_-}dz
%e^{-\frac{\vert z-(y-x-a_T^{-1} t \hat \mu)\vert^2}{2a_T^{-1} t}}
% \\
%&\qquad +e^{\langle y-x, \hat{\mu} \rangle}\int_{\mathbb{R}_-\times \mathbb{R}_+}dz
%e^{-\frac{\vert z-(y-x+a_T^{-1} t \hat \mu)\vert^2}{2a_T^{-1} t}}+e^{\langle y-x,\mu \rangle} \int_{\mathbb{R}_-^2}dz
%e^{-\frac{\vert z-(y-x+a_T^{-1}t \mu)\vert^2}{2a_T^{-1} t}}\bigg).
%\end{split}
%\label{eq:8}
%\end{equation}
Then for each of the 8 one-dimensional integral with respect to $z_i$, we proceed as in the one-dimensional case. That is, we split the integral into two cases depending whether the mean belongs to the same half line as $z_i$, $ i=1,2 $ or not. If the mean belongs to the  same interval as $z_i$, then  the integral can be lower bounded by $1/2$, otherwise, we apply the inequality $\vert a-b\vert^2\leq 2 (\vert a\vert^2+\vert b \vert^2)$ and use the fact that Gaussian integral with zero mean on the half line is equal to $1/2$.  Here, we will show how to deal with one of these terms, all other three follow similarly. In fact, for the second term of \eqref{eq:8}, we have 
\begin{align*}	&\frac{\mu_1 \mu_2}{8\pi tA_T }e^{\frac{a_T^{-1} t\vert  \mu \vert^2}{2 }}
e^{-\langle y-x, \hat{\mu} \rangle}
 \int_{\mathbb{R}_+\times \mathbb{R}_-}dz
e^{-\frac{\vert z-(y-x-a_T^{-1} t \hat \mu)\vert^2}{2a_T^{-1} t}}\\
&\qquad\qquad \geq \frac{\mu_1 \mu_2}{ 32 A_T a_T} e^{\frac{a_T^{-1} t\vert  \mu \vert^2}{2 }}e^{-\langle y-x, \hat{\mu} \rangle}
\left({\bf 1}_{y_1-x_1-\mu_1 t a_T^{-1}<0}  e^{-\frac{(y_1-x_1-\mu_1 t a_T^{-1})^2}{a_T^{-1} t}}+ {\bf 1}_{y_1-x_1-\mu_1 t a_T^{-1}\geq  0 }\right)\\
&\qquad\qquad \times
 \left({\bf 1}_{y_2-x_2+\mu_2 t a_T^{-1}>0}  e^{-\frac{(y_2-x_2+\mu_2 t a_T^{-1})^2}{a_T^{-1} t}}+ {\bf 1}_{y_2-x_2+\mu_2 t a_T^{-1}\leq  0 }\right).
\end{align*}
%This yields
%\begin{equation*} 
%\begin{split}
%&q_t(x,y)\\
%&\geq C_T e^{\frac{a_T^{-1} t\vert  \mu \vert^2}{2 }} \bigg\{e^{- \langle y-x, \mu \rangle} \left({\bf 1}_{y_1-x_1-\mu_1 t a_T^{-1}<0}  e^{-\frac{(y_1-x_1-\mu_1 t a_T^{-1})^2}{a_T^{-1} t}}+ {\bf 1}_{y_1-x_1-\mu_1 t a_T^{-1}\geq  0 }\right)\\
%&\qquad \times \left({\bf 1}_{y_2-x_2-\mu_2 t a_T^{-1}<0}  e^{-\frac{(y_2-x_2-\mu_2 t a_T^{-1})^2}{a_T^{-1} t}}+ {\bf 1}_{y_2-x_2-\mu_2 t a_T^{-1}\geq  0 }\right)\\
%&\qquad +
%e^{-\langle y-x, \hat{\mu} \rangle}
%\left({\bf 1}_{y_1-x_1-\mu_1 t a_T^{-1}<0}  e^{-\frac{(y_1-x_1-\mu_1 t a_T^{-1})^2}{a_T^{-1} t}}+ {\bf 1}_{y_1-x_1-\mu_1 t a_T^{-1}\geq  0 }\right)\\
%&\qquad \times \left({\bf 1}_{y_2-x_2+\mu_2 t a_T^{-1}>0}  e^{-\frac{(y_2-x_2+\mu_2 t a_T^{-1})^2}{a_T^{-1} t}}+ {\bf 1}_{y_2-x_2+\mu_2 t a_T^{-1}\leq  0 }\right)\\
%&\qquad +e^{\langle y-x, \hat{\mu} \rangle} \left({\bf 1}_{y_1-x_1+\mu_1 t a_T^{-1} > 0} e^{-\frac{(y_1-x_1+\mu_1 t a_T^{-1})^2}{a_T^{-1} t}}+ {\bf 1}_{y_1-x_1+\mu_1 t a_T^{-1} \leq  0 }\right)\\
%&\qquad \times \left({\bf 1}_{y_2-x_2-\mu_2 t a_T^{-1} < 0} e^{-\frac{(y_2-x_2-\mu_2 t a_T^{-1})^2}{a_T^{-1} t}}+ {\bf 1}_{y_2-x_2-\mu_2 t a_T^{-1}\geq  0 }\right)\\
%&\qquad +e^{\langle y-x,\mu \rangle}\left({\bf 1}_{y_1-x_1+\mu_1 t a_T^{-1} > 0}  e^{-\frac{(y_1-x_1+\mu_1 t a_T^{-1})^2}{a_T^{-1} t}}+ {\bf 1}_{y_1-x_1+\mu_1 t a_T^{-1}\leq   0 }\right)\\
%&\qquad \times \left({\bf 1}_{y_2-x_2+\mu_2 t a_T^{-1} >0}  e^{-\frac{(y_2-x_2+\mu_2 t a_T^{-1})^2}{a_T^{-1} t}}+ {\bf 1}_{y_2-x_2+\mu_2 t a_T^{-1}\leq   0 }\right) \bigg\},
%\end{split}
%\end{equation*}
From the above expression, it is easy to see that the above terms are estimated separately as in the one dimensional case. That is,  expanding the squares in the exponentials as in the one dimensional case, we conclude that   there exists $C_T>0$ such that
\begin{equation*} 
\begin{split}
q_t(x,y)\geq C_T e^{-\frac{a_T^{-1} t\vert  \mu \vert^2}{2 }} e^{-2\vert y_1-x_1\vert \mu_1}e^{-2\vert y_2-x_2\vert \mu_2}.
\end{split}
\end{equation*}

The rest of the proof for the case $d=2$ follows exactly as in the proof of the one-dimensional case. That is, appealing to formula (\ref{density}) yields for some positive  constants $\overline{C}_T$ and $c_T$,
	\begin{align*} 
f_t(x,y) &\geq \overline{C}_T e^{-(c_T+\lambda) t} e^{-4 \vert y_1-x_1 \vert \mu_1}e^{-4 \vert y_2-x_2 \vert \mu_2} \geq \overline{C}_T e^{-(c_T+\lambda) T} e^{-4 \vert y-x \vert \vert \mu \vert}.
\end{align*}

The same argument above can be extended to $\R^d$ for $d>2$ by splitting the integral into the $2^d$ orthants of $\R^d$. This completes the proof.
\end{proof}

Concerning  upper tail bounds, we observe that although Proposition \ref{p2} can be easily extended, Proposition \ref{p1} uses a one-dimensional argument which cannot be easily extended to the multidimensional setting. Thus, we leave it for further work.

\end{document}